\documentclass{amsart}                   
\usepackage{amsfonts}
\usepackage{amsmath}
\usepackage{amscd}
\usepackage{amssymb}
\usepackage[dvips]{graphicx}
\usepackage[dvips]{epsfig}
\usepackage{color}

\def\fraisse{Fra\"{\i}ss\'e\ }

\def\mf#1{\mathfrak{#1}}
\def\mb#1{\mathbb{#1}}
\def\mc#1{\mathcal{#1}}

\def\sub#1{_{_{#1}}}
\def\sup#1{^{^{#1}}}
\def\normal{T\sub4}

\def\Lim#1{\underset{#1}{\text{\rm Lim}}\ }

\def\R#1{\mb{R}\sup{#1}}                      
\def\s#1{\mb{S}\sup{#1}}                      

\def\im{\text{\rm Im}}              

\def\prof#1{\text{\rm len}(#1)}

\def\graph#1{\Gamma\sub{#1}}

\def\strat#1{\mc{S}\sub{#1}}           
\def\Strat#1{\left(#1,\strat{#1}\right)}

\def\ARROW#1{
\text{\begin{picture}(35,18)(15,15)         
            \put(30,22){$_{#1}$}
            \put(18,17){\vector(1,0){30}}
\end{picture}}}

\def\LARROW#1{
\text{\begin{picture}(35,18)(15,15)         
            \put(30,22){$_{#1}$}
            \put(48,17){\vector(-1,0){30}}
\end{picture}}}

\def\E{\mc{E}}

\def\Obj#1{\text{\rm Obj}\!\left(\mathfrak{#1}\right)}
\def\Mor#1{\text{\rm Hom}\!\left(\mathfrak{#1}\right)}

\def\Iso#1{\text{Iso}\!\left(\mathfrak{#1}\right)}

\newcounter{numero}

\newcounter{letra}
\newcommand{\Letra}{\medskip \setcounter{letra}{1}(\alph{letra}) }
\newcommand{\letra}{\medskip \addtocounter{letra}{1}(\alph{letra}) }

\newcounter{romnumero}

\newcounter{bibnumero}

\newtheorem{teo}{Theorem}[subsection]                  
\newtheorem{lema}[teo]{Lemma}
\newtheorem{prop}[teo]{Proposition}
\newtheorem{cor}[teo]{Corollary}
\newtheorem{quest}[teo]{Question}
\newtheorem{quests}[teo]{Questions}

\def\bteo{\begin{teo}}
\def\eteo{\end{teo}}
\def\bprop{\begin{prop}}
\def\eprop{\end{prop}}
\def\bcor{\begin{cor}}
\def\ecor{\end{cor}}
\def\blema{\begin{lema}}
\def\elema{\end{lema}}
\def\bquest{\begin{quest}}
\def\equest{\end{quest}}
\def\bquests{\begin{quests}}
\def\equests{\end{quests}}

\theoremstyle{definition}                           
\newtheorem{definition}[teo]{Definition}
\newtheorem{definitions}[teo]{Definitions}
\newtheorem{ejems}[teo]{Examples}
\newtheorem{ejem}[teo]{Example}

\def\bdeff{\begin{definition}\rm }
\def\edeff{\hfill$\square$ \end{definition}}
\def\bdefs{\begin{definitions}\rm }
\def\edefs{\end{definitions}}
\def\bejem{\begin{ejem}\rm }
\def\eejem{\end{ejem}}
\def\bejems{\begin{ejems}\rm }
\def\eejems{\end{ejems}}

\theoremstyle{remark}                               
\newtheorem{obs}[teo]{Remark}

\def\bobs{\begin{obs} }
\def\eobs{\end{obs}}

\newenvironment{dem}{ [{\it Proof\/}]\rm\hskip3mm }{\hfill$\square$\vskip5mm}       
\newenvironment{sketch}{ [{\it Sketch of the Proof\/}]\rm\hskip3mm }{\hfill$\square$\vskip5mm}

\def\bdem{\begin{dem}}
\def\edem{\end{dem}}

\def\bsketch{\begin{sketch}}
\def\esketch{\end{sketch}}

\begin{document}

\author{J. Mijares}
\address{Departamento de Matem\'aticas, Instituto Venezolano de Investigaciones Cient\'{\i}ficas.\\
Escuela de Matem\'aticas. Universidad Central de Venezuela.}
\email{jose.mijares@ciens.ucv.ve}

\author{G. Padilla}
\address{Departamento de Matem\'aticas, Universidad Nacional de Colombia, K30 con calle 45. 
Edificio 404, ofic. 315. Bogot\'a. (+571)3165000 ext 13166}
\email{gipadillal@unal.edu.co}

\title{On limit amalgamations of stratified spaces}

\date{10/04/2009}

\begin{abstract}
In this article we prove that stratified spaces and other geometric subfamilies  
satisfy categorical \fraisse properties, a matter that might be of interest for both geometers
and logicians. Part of this work was presented by the authors at the First Meeting of Logic and Geometry in Bogot\'a, on Sept. 2010.
\end{abstract}
\maketitle

\section*{Introduction}
All topological spaces considered in the sequel will be  Hlc2 (Hausdorff, locally compact 2nd countable) spaces. A stratified space is a topological space $X$ that can be decomposed in a family of disjoint pieces, the {\it strata}, which are locally closed manifolds satisfying nice incidence properties \cite{pflaum}.  Smooth manifolds are locally Euclidean spaces with smooth coordinate changes, and any smooth manifold has a trivial stratification whose singular part is the empty set.  In general, each point of a stratified pseudomanifold has a local neighborhood isomorphic to $\R{n}\times c(L)$ where $c(L)$ is the open cone of another compact stratified space, so an accurate definition is recursive. \vskip2mm

We are interested in Ramsey-like properties for stratified spaces and their limits. As a motivation, we show some examples of limits of stratified spaces which are finitely oscillation stable and so, in this approach, our first goal is to establish when a directed family of stratified spaces and embeddings has a limit. We also prove that the family of stratified pseudomanifolds is a \fraisse category.
There is a functor from stratified spaces to countable directed graphs. Basic stratified spaces correspond to finite connected graphs, and every stratified space can be obtained after a countable number of amalgamations and embeddings of basic ones. Under a suitable family of stratified embeddings we show that the family of stratified spaces is a \fraisse  category. The family of stratified pseudomanifolds also has this nice categoric behavior.\vskip2mm

\section{Bouquets of L\'evy families. A motivation}\label{bouquets}
Consider the sequence of unit spheres $\s{n}\cong\R{n}\sqcup\{\infty\}$ for $n\geq1$ whose singular part is $\Sigma\sub{n}=\{0,\infty\}$. 
The linear action of $U(n)$ on $\s{n}$ leaves $\Sigma\sub{n}$ fixed and is free elsewhere. The obvious equivariant smooth embeddings 
\[
		\s{n}\ARROW{}\s{n+1} 
\]
have as limit the unit sphere $\s{\infty}\subset l\sub2$, which is a smooth implicit Hilbert manifold \cite{lang}. The spheres $(\s{n},d\sub{n},\mu\sub{n})$, equipped with the Euclidean metric and the Haar measure, constitute a L\'evy family. This fact is one of the main ingredients in Milman's proof of the finite oscillation stability of the limit $\s{\infty}$, as a $U(l\sub2)$-space \cite{pestov}.\vskip2mm 

Let $X\sub{n}$ be a sequence of stratified pseudomanifolds and stratified morphism  $X\sub{n}\ARROW{\imath\sub{n}}X\sub{n+1}$. In general, we cannot guarrantee neither the existence of a decomposition of the limit $X\sub{\infty}=\Lim{n}X\sub{n}$ into disjoint smooth pieces nor the incidence properties of these pieces, even if the morphisms $\imath\sub{n}$ are stratified embeddings. In this article we find a class of strong embeddings, under which the limit space might be either
\begin{itemize}
		\item[$\bullet$] A decomposed space whose pieces are disjoint smooth (finite Euclidean or  Hilbert) manifolds. If the limit 
		length $\Lim{n}\prof{X}=\infty$ then the decomposition is not locally finite, 
		\item[$\bullet$] A Hilbert-smooth-stratified pseudomanifold if the graph of the stratification stabilizes beyond some $n$, 
		see \S\ref{def graph}. This is the case of the sphere $\s\infty$ above. The regular stratum is $l\sub2$, the link of the singular 	
		points $\{0,\infty\}$ is still $\s\infty$, which is not compact; or
		\item[$\bullet$] A stratified pseudomanifold, if $\Lim{n}\dim(X\sub{n})<\infty$ is finite.
\end{itemize}
Take a sequence 
\[
		G\sub1\subset \cdots \subset G\sub{n}\subset G\sub{n+1}\subset\cdots
\] 
of compact Lie groups. Suppose that each $X\sub{n}$ is a $G\sub{n}$-stratified pseudomanifolds (see \cite{gysin}) and each morphism $\imath\sub{n}$ is an equivariant strong embedding (see \S\ref{def special morphisms} below). By integrating with respect to the Haar measure of $G\sub{n}$ we can assume that each $X\sub{n}$ has a metric distance $d\sub{n}$ which is compatible with the smooth stratified structure and $G\sub{n}$ acts by stratified isometries \cite{bredon}. Choose, for each $n$, the probability measure $\mu\sub{n}$ induced by the normalized Haar measures of the group $G\sub{n}$ \cite{federer}. According to Theorem 1.2.10 in \cite{pestov}, if the spaces $(X\sub{n},d\sub{n},\mu\sub{n})$ constitute a L\'evy family, then  $X\sub{\infty}$ is a finite oscillation stable $G\sub{\infty}$-space,  where $G\sub{\infty}=\Lim{n}G\sub{n}$. Let us suppose that this is the case.\vskip2mm

Assume that for each $n$, the space $X\sub{n}$ has at least one fixed point, say $x\sub{n}$. Fix some integer $k\geq2$ and take the so called {\it bouquet} quotient space  
\[
Z\sub{n}=\frac{X\sub{n}\times[k]}{\{x\sub{n}\}\times\{1,\dots,k\}}
\]
where $[k]=\{1,\dots,k\}$ and the base point $z\sub{n}=\left[x\sub{n},j\right]$ is the class of $(x\sub{n},j)$. 
Take in $Z\sub{n}$ the unique distance $\widehat{d}\sub{n}$ that extends $d\sub{n}$ on each copy $X\sub{n}\sup{j}=q(X\sub{n}\times\{j\})$ where $q$ is the quotient map. Let $Z\sub{\infty}=\Lim{n}Z\sub{n}$ and consider the action of $G\sub{\infty}\sup{k}=G\sub{\infty}\oplus\overset{k\text{ times}}{\cdots}\oplus G\sub{\infty}$  on $Z\sub\infty$ given by 
\[
				(g\sub1,\dots,g\sub{k})([z,j])=[g\sub{j}(z),j]
\]
Fix some $\varepsilon>0$, an uniformly continuous function $Z\sub{\infty}\ARROW{f}\R{}$ and a finite subset $F\subset Z\sub{\infty}$. 
Define $F\sup{j}=\left(F\cup\left\{z\sub\infty\right\}\right)\cap X\sub\infty\sup{j}$. By the finite oscillation stability of $X\sub\infty$, for each $j\in[k]$ there is some $g\sub{j}\in G\sub\infty$ such that 
\[
		\text{diam}\left(g\sub{j}F\sup{j}\right)= 
		\text{sup}\left\{\ |f(q(x))-f(q(y))| : \ x,y\in g\sub{j}\left(F\sup{j}\right) \
		 \right\}<\varepsilon/2
\] 
We deduce that $\text{diam}\left(g F\right)<\varepsilon$ for $g=(g\sub1,\dots,
g\sub{k})$. This shows that $Z\sub\infty$ is finite oscillation stable under the stratified action of $G\sub\infty\sup{k}$. \vskip2mm

Notice that $\{Z\sub{n}\}\sub{n}$ is not a L\'evy family, since any probability measure in $Z\sub{n}$ comes from a convex combination
of the respective measures in the copies $X\sub{n}\sup{j},\ j=1,\dots,k$ of $X\sub{n}$. 

\section{\fraisse categories}
Usual \fraisse limits are formulated in terms of languages and models of finite structures. Given a certain language, a \fraisse family is, broadly speaking, a set of models and embeddings satisfying nice properties (heritability, joint embeddings, amalgamation and a countable number of non isomorphic models) which guarantees the existence of
a limit; see for instance \cite{chang,fraisse}. There is also a categoric treatment of these notions \cite{kubis}. In order to know
more about the connections of \fraisse theory with Ramsey theory and topological dynamics see \cite{kechris,pestov}.

\subsection{Embeddings and sub-objects}\label{def sub-objects}
Along this article we will work on a geometric category $\mf{C}$, i.e. 
a non-plain subcategory of $\mf{Top}$. The simple idea is that 
an {\bf embedding} is an arrow in a nice family of morphisms  $\E\supset\Iso{\mf{T}op}$; which 
satisfies some categoric properties.
Among them, for instance, embeddings are monomorphisms and the composition of embeddings is an embedding. 
We also say that $a$ is a {\bf sub-object} of $b$ if there is an embedding $a\ARROW{f}b$.
We will not stress on the categorical formalization of these notions, for more details see
\cite{gabriel,gelfand}.

\subsection{\fraisse categories}\label{fraisse categories}
	A diagram is said to be {\bf in} $\mf{C}$ if all its objects are in $\Obj{C}$ and all its arrows
	are in $\Mor{C}$.
		We say that $\mf{C}$ is a {\bf Fra\"iss\'e category} iff all its morphisms are embeddings and
	it satisfies the following axioms,
	\begin{enumerate}
		\item \underline{Heritability}: Each embedding 
			$a\ARROW{}b$ with target object $b\in\Obj{C}$ is a diagram in $\mf{C}$; 
			so $a\in\Obj{C}$.
		\item \underline{Joint embeddigs}: For any  $a,b\in\Obj{C}$ there is a diagram 
		\[
			a\ARROW{}c\LARROW{}b
		\]
		in $\mf{C}$. The object $c\in\Obj{C}$ is a {\bf joint object}.
		\item \underline{Amalgamation}: Each diagram $c\LARROW{}a\ARROW{}b$ in $\mf{C}$ 
		can be completed to  a commutative square 
		 \begin{center}
    			\begin{picture}(70,70)(15,15)        
            			\put(5,10){$c$}
            			\put(60,10){$d$}
            			\put(5,70){$a$}
            			\put(60,70){$b$}
 		             \put(7,65){{\vector(0,-1){40}}}
            			\put(62,65){{\vector(0,-1){40}}}
            			\put(15,75){\vector(1,0){40}}
            			\put(15,15){\vector(1,0){40}}
        \end{picture}
    		\end{center}\vskip2mm
		in $\mf{C}$. The object $d\in\Obj{C}$ is an {\bf amalgamation object}.
    \end{enumerate}

\subsection{Examples}\label{examples fraisse categories}
	 Here there are some examples of \fraisse categories.
	\begin{enumerate}
		\item The category of topological spaces and continuous functions $\mf{T}op$
		is a \fraisse category, since any two topological spaces can be topologically
		embedded in their disjoint union; and for any pair of topological embeddings
		$W\LARROW{h}X\ARROW{f}Y$ we can take the amalgamated sum
		\[
				Z=\frac{W\sqcup Y}{\sim}
				\hskip2cm
				f(x)\sim h(x)\forall x\in X
		\]
		The inclusions $Y\ARROW{}W\sqcup Y$ and $W\ARROW{}W\sqcup Y$ induce a commutative 		square of topological embeddings
		 \begin{center}
    			\begin{picture}(70,70)(15,15)  \label{diag1}       
            			\put(5,10){$W$}
            			\put(60,10){$Z$}
            			\put(5,70){$X$}
            			\put(60,70){$Y$}
 		             \put(7,65){{\vector(0,-1){40}}}
            			\put(62,65){{\vector(0,-1){40}}}
            			\put(15,75){\vector(1,0){40}}
            			\put(15,15){\vector(1,0){40}}
				    \put(0,45){$_h$}
      			    \put(66,45){$_\imath$}
	       		    \put(35,80){$_f$}
	        		    \put(35,7){$_\jmath$}
        		\end{picture}
    		\end{center}\vskip4mm
		What's more, $Z$ is a push-out since, for any other
		$Z'$ and any pair of continuous functions $W\ARROW{s}Z'\LARROW{t}Y$ inducing
		a commutative square diagram, i. e. such that $th=sf$; there is a {\it unique} continuous
		function $Z\ARROW{\phi}Z'$ such that $\phi\jmath=s$, $\phi\imath=t$; and if
		$s,t$ are topological embeddings then so is $\phi$. In the sequel, we will prefer the
		notation $Z=W\underset{_X}{\cup}Y$ for the amalgamated sum.
	\item The family of $CW$-complexes is closed under disjoint unions, embeddings and 
	topological amalgamations of (suitable) $CW$-embeddings. It is therefore a Fra\"iss\'e 
	category.   
	\item The family of smooth manifolds is not closed under under smooth amalgamations. 
	The 8-curve, which has a singular point, can be obtained as the amalgamation of two disjoint circles $Y,W$ through a
	distinguished point $X=\{p\}$.
\end{enumerate}
 
\section{Stratified spaces}
In order to recover any geometric smooth sense for amalgamations we must allow singular points. 
One way is to consider them as objects in a larger category, 
the family of smooth stratified spaces and its morphisms \cite{pflaum}. 

\subsection{Stratified spaces}\label{def stratified spaces}
    Let $X$ be a Hlc2 topological space.  A {\bf stratification} of $X$ is a locally finite partition $\strat{}$ 
    of $X$. The elements of $\strat{}$ are called {\bf strata}, and they are disjoint locally compact smooth manifolds 
     satisfying an 
    {\bf incidence condition}: Given any two strata $S,S'\in\strat{}$, if  $\overline{S}\cap S'\neq\emptyset$ then $S'
    \subset\overline{S}$ and we say that $S'$ {\bf adheres} to $S$ or just $S'\leq S$.
	If $\strat{}$ is a stratification of $X$,  we say that $(X,\strat{})$ is a {\bf stratified space}, although
	we might talk about {\it "the stratified space $X$"} when the choice of $\strat{}$ is clear in the context. 
Given a stratified space $(X,\strat{})$, 
\begin{enumerate}
	  \item The incidence condition is partial order relation.
	  \item For each stratum $S\in\strat{}$;\vskip2mm 
	  \begin{enumerate}
	  	\item $S$ is maximal (resp.  minimal) if and only if it is open (resp. closed).
  	  	\item The closure of $S$ is the union of the strata which adhere to it, 
		$\overline{S}=\underset{_{S'\leq S}}{\bigsqcup}\ S'$.
  	  	\item The set $U_S=\underset{_{S\leq S'}}{\bigsqcup}\ S'$ is open, we call it
          	the {\bf incidence neighborhood} of $S$.
	  \end{enumerate}
\end{enumerate}
\vskip2mm
Since $\strat{}$ is locally finite in $X$,  given a stratum $S$ 
the maximal length of all strict order chains in $\strat{}$ starting at $S$ 
is finite. The length of $S$ is the maximal integer $l=\prof{S}$ such that there is some strict 
order chain $S=S\sub0<\cdots<S\sub{l}$ in $\strat{}$. The {\bf length} of $X$
is the supremum of the lengths of the strata,
we write it $\prof{X}$. \vskip2mm
     
     A stratum $S$ is maximal (resp. minimal) iff it is open (resp. closed). Open strata will be called 
     {\bf regular} and, by opposition, a  {\bf singular} stratum will be a non open one.  
     The {\bf singular part}  $\Sigma$ (resp. {\bf minimal part} 
     $\mc{M}$) is the union of all singular (resp. minimal) strata. 
    \vskip2mm
	A {\bf morphism} $X$\ARROW{f}$Y$ between stratified spaces $\Strat{X}$ and $\Strat{Y}$ is a continuous function that sends smoothly each stratum of $X$ to some stratum of $Y$. The induced arrow	$\left(\strat{X},\leq\right)$\ARROW{f^{*}}$\left(\strat{Y}\leq\right)$ is a poset morphism.  An {\bf isomorphism} is a	bijective morphism whose inverse is also a morphism, the
$f^{*}$ induced by an isomorphism $f$ is a poset isomorphism.

\subsection{Examples of stratified spaces}\label{ejems stratified spaces}
    
	\begin{enumerate}
		\item Each manifold $M$ is a stratified space with respect to the family $\strat{M}=\{M\}$,
		we call it the trivial stratification  of $M$.
		\item The Cartesian product $X\times Y$ of stratified spaces $\Strat{X}$ and 
		$\Strat{Y}$ is a stratified space, with the canonical stratification
		\[
			\strat{X\times Y}=\left\{S\times T:S\in\strat{X},T\in\strat{Y}\right\}
		\]
		\item The disjoint union $X\sqcup Y$ of stratified spaces $\Strat{X}$ and 
		$\Strat{Y}$ is a stratified space, with the canonical stratification
		\[
			\strat{X\sqcup Y}=\left\{S:S\in\strat{X}\ \text{\rm or }S\in\strat{Y}\right\}
		\]
		\item Any open subset of a stratified space is also a stratified space.

		\item  Given a compact stratified space $\Strat{L}$, the {\bf open cone} $c(L)=\frac{L
		\times[0,\infty)}{L\times\{0\}}$ has a canonical stratification  
		\[
				\strat{c(L)}=\{v\}\cup\{S\times\R{+}:S\in\strat{L}\}
		\]
		We write $[p,r]$ for the equivalence class of a point $(p,r)$, and $v$ for the equivalence class 
		of $L\times\{0\}$, which we call the {\bf vertex} of the cone. We also adopt the convention 
		$c(\emptyset)=\{\star\}$.
\end{enumerate}

\subsection{Stratified subspaces}\label{def stratified subspaces}
    Let  $(X,\strat{})$ be a stratified space and $Z\subset X$ any topological subspace.
    The family 
    \[
    		\strat{Z}=\left\{S\cap Z:S\in \strat{}\right\}
	\]	
	 is the {\bf induced partition} of $Z$. We will also consider the family
	 \[
	 		\strat{(X,Z)}=\left\{S\cap Z,\ S\cap\left(\partial Z\right), \left(S-\overline{Z}\right)\ :\
			S\in\strat{} \right\}
	\]
	where $\partial Z=\overline{Z}\cap\overline{X-Z}$ 
	is the topological boundary. This family is the {\bf refinement} of 
	$\strat{X}$ induced by $Z$.  We will say that $Z$ is a {\bf stratified 
	subspace} (resp. a {\bf regular stratified subspace}) of $X$ iff $\strat{Z}$ is a stratification 
	of $Z$ (resp. iff $\strat{(X,Z)}$ is a stratification of $X$), which happens iff the intersection of 
	$Z$ (resp. and $\partial Z$) with any stratum $S\in\strat{}$ is a submanifold of $S$.

\subsection{Embeddings}\label{def special morphisms}
    Let $X\ARROW{f} Y$ be a morphism; we say that...
    \begin{enumerate}
			 \item $f$ is an {\bf immersion} 
			iff the restriction $f\mid_{_S}$ to each stratum
			 $S\in\strat{X}$ a smooth immersion. 
    			\item $f$ is a {\bf weak embedding} iff it is a 1-1 immersion satisfying the {\bf lifting 
			property}:  For each morphism $X'\ARROW{h} Y$ such that $\im(h)\subset\im(f)$, the 			composition $\widehat{f}=fh\sup{-1}$ is a morphism $X'\ARROW{\widehat{f}}X$. 
    			\item $f$ is an {\bf embedding} iff $f(X)$ is a stratified subspace of $Y$ and 
			$X\ARROW{f} f(X)$ is an isomorphism. 
			\item  $f$ is a {\bf strong embedding} iff $f$ is an embedding, $f(X)$ is regular and 
			$\strat{Y}=\strat{(Y,f(X))}$.
    \end{enumerate}

\subsection{Examples}\label{example morphisms}
		 \hfill
	\begin{enumerate}
		 \item Each embedding is a weak embedding.
		 \item If $X\ARROW{f}Y$ is a weak embedding, then 
		 the stratification of $X$ is uniquely determined by $Y$ up to isomorphisms. 
		\item $Z\subset X$ is a stratified subspace (resp. regular) iff the inclusion 
		$Z\ARROW{}X$ is an embedding (resp. a strong embedding). 
		\item Consider the 8-curve 
		$\gamma\subset\R2$ and let $p\in\gamma$ be the singular 
		point. Let $C\sub1,C\sub2$ be the two connected components of $\gamma-p$.
		Consider in $\gamma$ the stratifications $\strat0=\{p,(\gamma-p)\}$
		and $\strat1=\{p,C\sub1,C\sub2\}$. Also consider in $\R2$
		the stratifications $\mc{T}\sub0=\{\R{2}\}$, $\mc{T}\sub1=\{p,(\R{2}-p)\}$
		$\mc{T}\sub2=\{p,(\gamma-p),(\R{2}-\gamma)\}$ and 
		$\mc{T}\sub3=\{p,C\sub1,C\sub2,(\R{2}-\gamma)\}$.
		 Write $\gamma\sub{i}$ and
		$\R2\sub{j}$ for the stratified spaces $(\gamma,\strat{i})$, and $(\R2,\mc{T}\sub{j})$. 
		Then\vskip2mm
		\begin{itemize}
			\item The identity $\gamma\sub{k}=\gamma\sub{l}$ is not a morphism for $k=0,l=1$; 
			an embedding for $k=1,l=0$ and $k=0$; and an isomorphism for $k=l$.\vskip2mm
			\item The inclusion $\gamma\sub{k}\subset\R2\sub{l}$ is a proper (1-1)-immersion for $k\leq l$.
			It is an embedding if $k=l+1$, and a strong embedding if $k=l+2$. For
			$k=0$ and $l=3$ it is not a morphism.\vskip2mm
		\end{itemize}
		A more precise explanation of this situation is given below.
\end{enumerate}

\subsection{Reducibility}\label{def reducibility}
	The disjoint union of two stratified spaces is again a stratified space, see 
	\S\ref{ejems stratified spaces}-(3). 
	A stratified space $(X,\strat{})$ is {\bf reducible} if it can be written as the disjoint union
	of two stratified spaces; and {\bf irreducible} if not. If $X=X\sub1\sqcup X\sub2$ is such a disjoint union
	then, by \S\ref{def stratified subspaces} $X\sub1,X\sub2$ are regular stratified subspaces of $X$.\vskip2mm

	An {\bf irreducible component} of $X$	is a minimal irreducible (regular) stratified subspace
	of $X$. If the strata of $X$ are all connected manifolds, then an irreducible component is a connected
	component of $X$; although in general irreducible components might not be connected. 
 
\subsection{Associated graphs}\label{def graph} 
	Given a stratified space $(X,\strat{})$;
	the {\bf graph} $\graph{X}$ \linebreak associated to $X$ is the 
	directed graph induced by the poset $(\strat{},\leq)$.\vskip2mm 
	A {\bf vertex} of $\graph{X}$ is a stratum $S\in\strat{}$. A {\bf directed edge} $\{S\sub2,S\sub1\}$ 
	is a minimal strict 
	incidence chain	$S\sub1<S\sub2$, which means that there are no intermediate strata between $S\sub1,S\sub2$.
 
\hrulefill\vskip1mm
    \begin{center}
	
    \begin{picture}(100,110)(15,15)  \label{diag2}       
            \put(-100,100){$\R2\sub0$}
            \put(-80,100){$\bullet$}

            \put(-40,100){$\R2\sub1,\gamma\sub0$}
            \put(-2,110){$\bullet$}
            \put(0,110){{\line(0,-1){20}}}
            \put(0,110){{\vector(0,-1){10}}}
            \put(-2,110){$\bullet$}
            \put(-2,90){$\bullet$}

            \put(60,100){$\R2\sub2$}
            \put(78,120){$\bullet$}
            \put(80,120){{\line(0,-1){20}}}
            \put(80,120){{\vector(0,-1){10}}}
            \put(78,100){$\bullet$}
            \put(80,100){{\line(0,-1){20}}}
            \put(80,100){{\vector(0,-1){10}}}
            \put(78,80){$\bullet$}

            \put(140,100){$\gamma\sub1$}
            \put(165,109){$\bullet$}
            \put(212,109){$\bullet$}
            \put(170,110){{\line(1,-1){20}}}
            \put(170,110){{\vector(1,-1){10}}}
            \put(212,110){{\line(-1,-1){20}}}
            \put(212,110){{\vector(-1,-1){10}}}
            \put(189,87){$\bullet$}

            \put(-70,30){$\R2\sub1$}
            \put(-22,51){$\bullet$}
            \put(-17,52){{\line(1,-1){20}}}
            \put(-17,52){{\vector(1,-1){10}}}
            \put(-21,52){{\line(-1,-1){20}}}
            \put(-21,52){{\vector(-1,-1){10}}}
            \put(-45,29){$\bullet$}
            \put(1,29){$\bullet$}
            \put(-40,30){{\line(1,-1){20}}}
            \put(-40,30){{\vector(1,-1){10}}}
            \put(2,30){{\line(-1,-1){20}}}
            \put(2,30){{\vector(-1,-1){10}}}
            \put(-21,7){$\bullet$}

            \put(30,35){$_{[0,1]}$}
            \put(78,51){$\bullet$}
            \put(82,52){{\line(1,-1){20}}}
            \put(82,52){{\vector(1,-1){10}}}
            \put(78,52){{\line(-1,-1){20}}}
            \put(78,52){{\vector(-1,-1){10}}}
            \put(55,29){$\bullet$}
            \put(101,29){$\bullet$}

            \put(130,35){$_{[0,1]\sup{2}}$}
            \put(160,0){$\bullet$}
            \put(160,30){$\bullet$}
            \put(160,60){$\bullet$}
            \put(220,0){$\bullet$}
            \put(220,30){$\bullet$}
            \put(220,60){$\bullet$}
            \put(190,0){$\bullet$}
            \put(190,30){$\bullet$}
            \put(190,60){$\bullet$}
            \put(163,1){{\line(0,1){60}}}
            \put(193,1){{\line(0,1){60}}}
            \put(223,1){{\line(0,1){60}}}
            \put(161,2){{\line(1,0){60}}}
            \put(161,32){{\line(1,0){60}}}
            \put(161,62){{\line(1,0){60}}}
            \put(193,32){{\vector(0,1){10}}}
            \put(193,32){{\vector(0,-1){10}}}
            \put(193,32){{\vector(1,0){10}}}
            \put(193,32){{\vector(-1,0){10}}}
            \put(163,32){{\vector(0,1){10}}}
            \put(163,32){{\vector(0,-1){10}}}
            \put(223,32){{\vector(0,1){10}}}
            \put(223,32){{\vector(0,-1){10}}}

            \put(193,2){{\vector(1,0){10}}}
            \put(193,2){{\vector(-1,0){10}}}
            \put(193,62){{\vector(1,0){10}}}
            \put(193,62){{\vector(-1,0){10}}}
        \end{picture}\vskip5mm\hrulefill\vskip2mm
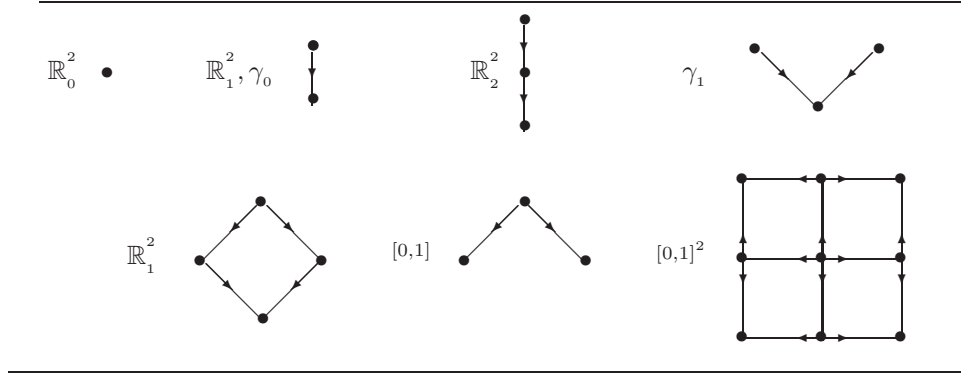
	\tiny{{\bf Figure 1:} Some examples of graphs induced by stratifications, cf. \S\ref{example morphisms}-(4)}
    \end{center}

\subsection{Properties of the associated graph}\label{properties of graphs}
Given a stratified space $(X,\strat{})$.
\begin{enumerate}
		\item $\graph{X}$ has (at most) a countable number of vertices. 
		\item Each path in $\graph{X}$ is finite. 
		\item $\prof{X}$ is the (supremum) of the lengths of the paths in $\graph{X}$. 
		\item For each stratum $S\in\strat{}$
		\begin{enumerate}
			\item $\graph{\overline{S}}$ is the subgraph in $\graph{X}$ consisting 
			of all paths starting at $S$.
			\item $\graph{U\sub{S}}$ is the finite subgraph in $\graph{X}$ consisting 
			of all paths ending at $S$.
		\end{enumerate}
		\item $X$ is irreducible if and only if $\graph{X}$ is a connected graph.
		\item The graph of a normal stratified pseudomanifold is a tree 
		(\footnote{Stratified normality is the property that
		all the links of $X$ are connected; it should not be confused with 
		the $\normal$ separation axiom, which is satisfied by 
		any stratified space. See the next sections and \cite{normalizer}
		for more details. }).
		\item 	A regular subspace $Z\subset X$ is closed if and only if $\graph{Z}$ is a 
		subgraph of $\graph{X}$. Strong embeddings of stratified spaces correspond to
		strong embeddings of graphs in the sense of \cite[3.4-p.414]{hrushovski}, so the family has
		a partial order that takes into account the embeddings of 
		regular subspaces \cite[p.364]{fraisse}.
\end{enumerate}
These properties can be easily checked, we leave the details to the reader.

\subsection{Amalgamations of stratified spaces}\label{subsection amalgamations of stratified spaces}
Stratified spaces and strong embeddings provide a full geometric sense of smoothness 
for amalgamations.
\blema\label{lema proper 1-1 immersions}
	Let $X\ARROW{f}Y$ be a proper 1-1 immersion. Then $f$ is an embedding
	iff the induced poset morphism
	$\left(\strat{X},\leq\right)\ARROW{f^{*}}\left(\strat{Y},\leq\right)$
	 is 1-1.
\elema
\bdem
	Notice that
	\begin{enumerate}
		\item \underline{$X\ARROW{f}f(X)$ is an homeomorphism:} This is a consequence of two facts,
			\begin{enumerate}
					\item $X$ is locally compact and $Y$ is Hausdorff.
					\item A proper continuous bijection from a compact space to a Hausdorff spaces is a 
					homeomorphism.
			\end{enumerate}
		\item \underline{For each stratum $S\in\strat{X}$ the restriction $S\ARROW{f}f(S)$ 
		is a diffeomorphism:} Let $R\in\strat{Y}$ be the corresponding stratum such that 
		$f(S)\subset R$.
		Since $f\mid_{_S}$ is a smooth proper 1-1 immersion, we deduce that $S\ARROW{f}R$ is
		a smooth embedding, so $f(S)$ is a regular submanifold of $R$. Now $S\ARROW{f}f(S)$
		is a proper embedding between equidimensional manifolds, so it is a diffeomorphism.
		\item \underline{$f(X)$ is a stratified subspace of $Y$:} This is a consequence of the previous 
		step. The family $\strat{f(X)}$ is a stratification of $f(X)$.

		\item \underline{$X\ARROW{f}f(X)$ is an isomorphism iff $f^{*}$ is 1-1:} Write $h=f^{-1}$. 
		Since $f$ is a homeomorphism on its image, 
		the inverse map $f(X)\ARROW{h}X$ is continuous. The map $h$ sends each stratum of $f(X)$
		in a disjoint union of non comparable strata in $X$. What's more, by step (2), 
		the restriction of $h$ to each
		stratum of $f(X)$ is smooth.  Finally, $h$ is stratum preserving iff each stratum
		of $Y$ meets $f(X)$ in only one stratum, i. e. iff the induced map 
		$\Strat{f(X)}\ARROW{h}\Strat{X}$ is well defined, which happens iff $f^{*}$ is 1-1. 
	\end{enumerate}
\edem
As a consequence,
\bteo\label{lema stratified spaces are Fraisse}\hfill
	\begin{enumerate}
		\item The amalgamation of two stratified spaces by a pair of strong embeddings is 
		a stratified space.
		\item Stratified spaces constitute a \fraisse category with respect to strong embeddings.
	\end{enumerate}
\eteo
\bdem
		Statement (2) is a direct consequence of (1), which we now prove.\vskip2mm

		Any two stratified spaces $(W,\strat{W})$ and $(Y,\strat{Y})$ can be strongly 
		embedded in their disjoint union, so we have joint embeddings, see 
		axiom \S\ref{fraisse categories}-(1). We now verify \S\ref{fraisse categories}-(2),
		the amalgamation property.  Let $(X,\strat{X})$ be any other stratified space and assume 
		that the arrows $f,h$ in \S\ref{examples fraisse categories}-(1) are strong embeddings. Then
		\begin{itemize}
				\item \underline{The amalgamated sum $W\underset{_X}{\cup}Y$  is a stratified space:} 
				Let  $W\sqcup Y\ARROW{q}Z$ be the quotient map. The family
			   \[
					\strat{W\underset{_X}{\cup}Y}=\left\{q(S):S\in\strat{W}\ \text{\rm or }S\in
						\strat{Y}\right\}
			  \]
			     is a locally finite partition. Since $f,h$ are strong embeddings,  $X$ can be seen
			     as a regular stratified subspace of $W$ and $Y$  an the same time, and therefore
			     the elements of $\strat{W\underset{_X}{\cup}Y}$ are locally closed manifolds (with
			     the induced topology) and satisfy the incidence condition \S\ref{def stratified spaces}.
			     \item \underline{The quotient map $W\sqcup Y\ARROW{q}Z$ is a morphism:}
			     This is straightforward.
			     \item \underline{\S\ref{examples fraisse categories}-(1) is
			     a diagram of stratified  strong embeddings:} 
			     Since the inclusions of $W,Y$ in the disjoint union $W\sqcup Y$
			    are strong embeddings, the obvious induced maps 
			    \[
			    			W\ARROW{\jmath}W\underset{_X}{\cup}Y\LARROW{\imath}Y
				\] 
				are strong embeddings.	    
			\end{itemize}
			 This concludes the proof.
\edem
\bcor\label{cor amalgamation is minimal}
	The amalgamated sum of stratified spaces with strong embeddings is a pushout.
\ecor
This is an easy consequence, we leave the details to the reader.

\subsection{Basic stratified spaces}\label{def basic spaces}
 	A stratified space $(X,\strat{})$ is  {\bf basic} iff
	\begin{enumerate}
		\item It is irreducible.
		\item It has finite length.
	\end{enumerate}
	Therefore, $X$ is basic iff $\graph{X}$ is a connected finite graph.
 
\blema\label{lema amalgamation of basics} 
	Let $W\underset{_X}{\cup}Y$ be the amalgamation of two stratified spaces 
	along a closed regular subspace $X$. 
	\begin{enumerate}
		\item If $Y,W$ are irreducible, then so is $W\underset{_X}{\cup}Y$.
		\item If $Y,W$ are basic, then so is $W\underset{_X}{\cup}Y$.
	\end{enumerate}	
\elema
\bdem
	Let $W\LARROW{f}X\ARROW{h}Y$ be two strong embeddings, and assume that
	$X$ is a closed subspace of (both) $W,Y$. Then 
	the graph $\graph{W\underset{_X}{\cup}Y}$ of the amalgamated space is the joint of 
	$\graph{W}$ and $\graph{Y}$  along $\graph{X}$; i. e. 
	\[	
		\graph{W\underset{_X}{\cup}Y}=\graph{W}\underset{_{\graph{X}}}\vee\graph{Y}
	\]
	We conclude that if $\graph{W},\graph{Y}$ are connected (resp. finite) graphs then so
	is $\graph{W\underset{_X}{\cup}Y}$. 
\edem

Here there is another useful and easy result.
\blema\label{lema noncomparable strata}
	Non comparable strata can be separated with disjoint
	open subsets.
\elema
\bdem
	Notice that\vskip2mm
	\begin{itemize}
		\item[\Letra] \underline{It is enough to show it for minimal strata:}
		If $\mc{F}\subset\strat{}$ is a family of non-comparable strata, take
		the union of the incidence neighborhoods
		$Z=\underset{_{S\in\mc{F}}}{\cup} U\sub{S}$. Then $Z$ is open in $X$,
		and $S\in\mc{F}$ iff $S$ is a minimal stratum in $Z$. 
		Since $Z$ is open, it is enough to give a family of disjoint neighborhoods
		in $Z$ separating the strata in $\mc{F}$.
		\item[\letra] \underline{ Minimal strata can be separated by disjoint open subsets}: 
		Any two different minimal strata in $X$ are disjoint closed subsets, that can be separated
		with two disjoint open subsets because $X$ is $\normal$.	
		The whole family of minimal strata can be separated because of \S\ref{def stratified spaces}-(1), 
		(2)-(b), and the facts that
		$X$ is $\normal$, and $\strat{}$ is locally finite \cite{tesis}.
	\end{itemize}
\edem

\bprop\label{prop limits of stratified spaces}
	Each stratified space is the result of, at most, a countable number of 
	disjoint unions or amalgamations of basic ones.
\eprop
\bdem
	Let $(X,\strat{})$ be a stratified space. 
	Take a minimal stratum $S$ in $X$. By\S\ref{properties of graphs}-(4b) and  
	\S\ref{def basic spaces}; $U\sub{S}$ is basic.
	If $S'\neq S$ is another minimal stratum and $U\sub{S}\cap U\sub{S'}\neq\emptyset$
	is non-empty; then the graph $\graph{U\sub{S}\cup U\sub{S'}}$ 
	is connected and 
	$U\sub{S}\cup U\sub{S'}$ is irreducible.   
	(\footnote{Notice that for any stratum $T$ such that $T\subset U\sub{S}\cup U\sub{S'}$
	we have $U\sub{T}\subset U\sub{S}\cup U\sub{S'}$
	so the family $\{U\sub{S}:S\in\strat{}\}$ is a basis.})	Since $\graph{U\sub{S}},\graph{U\sub{S'}}$
	are finite then so is $\graph{U\sub{S}\cup U\sub{S'}}$. By \S\ref{def basic spaces} 	
	we deduce that	
	\[
		U\sub{S}\cup U\sub{S'}\cong U\sub{S}\underset{_{U\sub{S}\cap U\sub{S'}}}{\cup} U\sub{S'}
	\]
	is basic. Since the stratification is locally finite, by \S\ref{lema noncomparable strata} and 
	\S\ref{def stratified spaces}-(2) there is at most a countable number of minimal strata. Since
	\[
		X=\cup\{U\sub{S}:S\text{ is minimal } \}	
	\]
	we conclude that $X$ is the result of, at most, a countable number of unions or amalgamations 
	of basic stratified spaces.
\edem

\section{Stratified Pseudomanifolds}

A stratified pseudomanifold is a stratified space together a family
of conic charts which reflect the way in which we approach
the singular part. The definition is given by induction 
on the length. 

\subsection{Stratified pseudomanifolds}\label{def spm} 
A 0-length {\bf stratified pseudomanifold} is a smooth manifold with the
trivial stratification. A stratified space $\left(X,\strat{}\right)$ with $l(X)>0$ 
is a {\bf stratified pseudomanifold} if, for each singular stratum $S$, 
\begin{enumerate}
	\item   There is a compact stratified pseudomanifold $\Strat{L}$ with $l(L)<l(X)$. We call 
	$L$ the {\bf link} of $S$ because
	\item   Each point $x\in S$ has an open neighborhood $x\in U\subset S$
	and a stratified embedding $U\times c(L)\ARROW{\alpha}X$ 
	on an open neighborhood of $x\in X$. 
\end{enumerate}
The image of $\alpha$ is called a {\bf basic neighborhood} of $x$. 
Notice that $\Im(\alpha)\cap S=U$. Without loss of generality, 
we assume that $\alpha(u,v)=u$ for each $u\in U$ (where $v$ is the vertex of $c(L)$, see 
\S\ref{ejems stratified spaces}-(5)). We summarize the above situation by saying that the pair
$(U,\alpha)$ is a {\bf chart} of $x$. The family of charts is an {\bf atlas} of $(X,\strat{})$.

\subsection{Examples} \label{ejm examples of stratified pseudomanifolds}\hfill
\begin{enumerate}
  \item A basic model $U\times c(L)$  is a stratified pseudomanifold if $\Strat{L}$ is a compact
  stratified pseudomanifold.
  \item If $M$ is a manifold and $X$  is a stratified pseudomanifold then
  $M\times X$ is a stratified pseudomanifold. 
  \item Every open subset of a stratified pseudomanifold is again a stratified pseudomanifold. 
  \item Since algebraic manifolds satisfy the Withney's conditions, every algebraic manifold is a stratified
   pseudomanifold \cite{pflaum}. 
\end{enumerate}

Now we will extend some results of the previous section.

\bprop\label{lema amalgamation of pseudomanifolds} 
	Let $W,X,Y$ be stratified pseudomanifolds with finite length. If
	\begin{enumerate}
		\item $W\LARROW{f}X\ARROW{h}Y$ are strong embeddings, and
		\item $X$ is closed in both $W,Y$;
	\end{enumerate}
	then $W\underset{_X}{\cup}Y$ is a stratified pseudomanifold.
\eprop
\bdem
	The amalgamated space $Z=W\underset{_X}{\cup}Y$ is a finite length
	stratified space. We check the existence of conical charts 
	\S\ref{def spm} for each $z\in Z$. Since $X$ is closed in $Y$; if 
	If $z\in (Z-W)\cong (Y-X)$, then by example 
	\S\ref{ejm examples of stratified pseudomanifolds}-(3)
	there is nothing to prove. The same holds for $z\in (Z-Y)$ so we must
	check it only for $z\in W\cap Y=X$.\vskip2mm

	Let us remark that this is a local problem. Since $z$ has conical charts in $W,X,Y$,
	and up to some minor details we assume that $W=S\times c(L)$, $Y=S\times c(L')$ and $X=S\times c(N)$
	are basic neighborhoods and $f(u,z)=h(u,z)=u$ for any $u\in S$.
	\vskip2mm
	We now proceed by induction on $l=\prof{X}$.
	\begin{itemize}
		\item \underline{Case $l=0$:} Then $N=\emptyset$ and $X=S$. The amalgamated sum
		\[
			W\underset{_X}{\cup}Y=[S\times c(L)]\underset{_S}{\cup}[S\times c(L')]
			= S\times c(L\sqcup N)
		\]
		is a basic neighborhood of $z\in S$, and the link of $z$ is the disjoint union
		of the links in $W,Y$.
		\item \underline{Inductive hypothesis:} We assume \S\ref{lema amalgamation of pseudomanifolds} 
		for any	triple $(W',X',Y')$ such that $\prof{X'}\leq l-1$.
		\item \underline{General case $l>0$:} Notice that, by the hypotheses of 
		\S\ref{lema amalgamation of pseudomanifolds}, the (image of the) pseudomanifold $N$ is closed
		in $L,L'$. Since $\prof{N}<\prof{X}$; we can apply the inductive hypothesis to $(L,N,L')$
		in order to get a stratified pseudomanifold  
		$L\underset{_N}{\cup}L'$. 
		Therefore
		\[
				W\underset{_X}{\cup}Y=S\times c\left(L\underset{_N}{\cup}L'\right)
		\]
		is a stratified pseudomanifold.
	\end{itemize}
\edem

\bobs\label{remark on closed amalgamations}
	The condition that $X$ is closed, in \S\ref{lema amalgamation of pseudomanifolds}-(2), cannot be
	avoided. As a simple counter-example, consider $W=[0,1]\times[0,1]$ a closed square, 
	$Y=[1,2)\times(0,1)$ and $X=\{1\}\times(0,1)$. There is a link at $p=(1,0)\in W\underset{_X}{\cup}Y$
	but it is not compact.
\eobs

\bprop\label{prop limits of pseudomanifolds}
	Each stratified pseudomanifold is the result of, at most, a countable number of 
	disjoint unions or amalgamations of basic pseudomanifolds. 
\eprop
\bdem
	By \S\ref{prop limits of stratified spaces} 
	$X$ is the (at most countable) union of the incidence neighborhoods $U\sub{S}$ of all its (minimal) 
	strata.	By \S\ref{ejm examples of stratified pseudomanifolds} these open neighborhoods are basic 
	pseudomanifolds.
\edem

We summarize the results of this section in the following

\bteo\label{teo psv are Fraisse}
		The family of stratified pseudomanifolds, with respect to strong embeddings, 
		is a \fraisse category. Moreover, it is
		the closure by amalgamations of the family of basic stratified pseudomanifolds.
\eteo

	\section*{Acknowledgments}
We would like to thank A. Villaveces for some helpful conversations, and
the referees of the journal Revista Colombiana de Matem\'aticas for their comments which have 
substantially improved the final version of this manuscript. We also thank the Departaments of Mathematics at the Instituto Venezolano de Investigaciones Cient\'ificas (IVIC) and the Universidad Nacional de Colombia (UNAL), for their support and hospitality.

\end{document}